\documentclass{article}
\usepackage{amssymb,amsmath}
\newtheorem{twr}{Theorem}
\newtheorem{lem}{Lemma}
\newtheorem{df}{Definition}
\begin{document}
\title{{\bf A Generalization of the Markov-Kakutani Fixed
Point Theorem}\footnote{The paper was supported by the 
KBN grant no. 5 P03B 09320}}
\author{Jaros\l{}aw Wawrzycki\footnote{Email address: 
Jaroslaw.Wawrzycki@ifj.edu.pl}\\ 
\footnotesize{Institute of Nuclear Physics PAS, ul. Radzikowskiego 152,}\\ 
\footnotesize{31-342 Krak\'ow, Poland}}
\date{\footnotesize{June 20, 2003}}
\maketitle
\begin{abstract} 
In this announcement we generalize the 
Markov-Kakutani fixed point theorem for abelian semi-groups 
of affine transformations extending it on some class of 
non-commutative semi-groups. As an interesting example we apply it 
obtaining a generalization of the invariant version of the 
Hahn-Banach theorem.
\end{abstract}

\begin{center}
\end{center}

\vspace{2ex}

We will consider a general linear topological space $X$,
a convex compact subset $K$ of $X$ and  semi-groups
of affine continuous transformations of $K$ into itself.
Further on we call them simply semi-groups.  
Specifically we confine our attention to a special family
$\mathcal{F}$ of such semi-groups.  

\begin{df} Let $\mathcal{F}$ be the smallest family of semi-groups
of affine transformations of $K$ fulfilling \\
(i) Every abelian group of affine transformations of $K$
belongs to $\mathcal{F}$.\\
(ii) If the semi-group $F$ possesses a normal factor $H$ and
$H$ together with $F/H$ belong to $\mathcal{F}$, then $F$ is an
element of $\mathcal{F}$ as well.\label{defin} \end{df}
We will show that each semi-group $F$ of $\mathcal{F}$ possesses
a fixed point $p$ in $K$: $p=f(p)$ for all $f \in F$.
We divide $\mathcal{F}$ into a sequence
$\mathcal{F}_{1}, \mathcal{F}_{2}, \ldots $ of subfamilies.
Semi-groups $F_{1}$ belonging to $\mathcal{F}_{1}$ are defined
to be abelian. Each semi-group $F_{2}\in \mathcal{F}_{2}$ 
is a semi-direct product of elements of $\mathcal{F}_{1}$, 
and in general each $F_{k} \in \mathcal{F}_{k}$ is a semi-direct 
product of semi-groups $F_{i}$ and $F_{j}$ 
with $i,j <k$. It will be sufficient to prove
our assertion for each class $\mathcal{F}_{k}$ of semi-groups 
separately.  

We need a number of Lemmas.
\begin{lem} If the semi-group $F$ is abelian, then\footnote{Sign 
$co(\centerdot)$ denotes the convex hull.} $co(F)$ is an abelian 
semi-group as well, and the family of sets
\begin{equation}\label{family}
\{f(K), f \in co(F)\}
\end{equation}
has the finite intersection property (any finite number of its 
elements has nonempty intersection).\label{lem1} \end{lem} 
PROOF. Firs of all $co(F)$ is a semi-group because\footnote{The
symbol $f\circ g$ denotes the composition of affine transformations 
$f$ and $g$.} for $f,g \in co(F)$ 
\begin{displaymath}
f\circ g = (\sum_{i} \alpha_{i}f_{i})\circ (\sum_{j} \beta_{j}g_{j})
= \sum_{i,j} \alpha_{i}\beta_{j}f_{i}\circ g_{j},
\end{displaymath}
and $f_{i}\circ g_{j}$ being element of $F$ as well as 
$\sum_{i,j} \alpha_{i}\beta_{j} = \sum_{i}\alpha_{i} 
= \sum_{j}\beta_{j} = 1$. Thus $f \circ g \in co(F)$. 
Obviously, $f\circ g = g \circ f$ transforms $K$ into $K$,
therefore $co(F)$ is abelian as well.
Let us take any two sets $f(K)$ and $g(K)$ of the family 
(\ref{family}). Then, $f(g(K)) \subset f(K)$ as well as 
$g(f(K)) \subset g(K)$, therefore
\begin{displaymath}
\emptyset \neq g\circ f(K) = f\circ g(K) \subset f(K) \cap g(K), 
\end{displaymath}
and a simple induction shows that the family (\ref{family})
has the finite intersection property. $\square$

\begin{lem} If $H$ is a normal factor of $F$, then
\begin{equation}\label{eqlem2}
co(H) \circ g=g\circ co(H), 
\end{equation}
for $g \in co(F/H)$.\label{lem2} \end{lem}
PROOF. We have $co(H) \circ g \ni f = (\sum_{i}\alpha_{i}h_{i})
\circ (\sum_{j}\beta_{j}g_{j}) 
= \sum_{i,j}\alpha_{i}\beta_{j}h_{i}\circ g_{j}$ with the 
appropriate $h_{i} \in H$ and $g_{j} \in F/H$. But $H$ is 
normal factor, so that $h_{i}\circ g_{j}=g_{j}\circ h'_{i}$
for some $h'_{i}\in H$. Therefore, 
$f= \sum_{j,i} \beta_{j}\alpha_{i}g_{j}\circ h'_{i}
= (\sum{j}\beta_{j}g_{j})\circ (\sum_{i}\alpha_{i}h'_{i}) 
\in g \circ co(H)$. $\square$ 
 
\begin{lem} Suppose $H$ to be normal factor of semi-group $F$,
and both $H$ and $F/H$ are abelian. Then, the family 
of sets 
\begin{equation}\label{family'}
\{f(K), f \in co(H)\circ co(F/H)\},
\end{equation}
has the finite intersection property.\label{lem3} \end{lem} 
PROOF. Note, please, that $co(H)\circ co(F/H)$ is a semi-group:
for $f_{1} = h_{1}\circ g_{1}$ and $f_{2} = h_{2}\circ g_{2}$ 
lying in  $co(H)\circ co(F/H)$ we have  $f_{1} \circ f_{2}
= h_{1}\circ g_{1}\circ h_{2} \circ g_{2}$. This is equal to
$h_{1}\circ h'_{2}\circ g_{1}\circ g_{2}$, for some 
$h'_{2}\in co(H)$, on account of Lemma \ref{lem2}. By Lemma
\ref{lem1} $h_{1}\circ h'_{2} \in co(H)$ and   
$g_{1}\circ g_{2} \in co(F/H)$, which proves that 
$co(H)\circ co(F/H)$ is a semi-group. Similarly, for
each $h\in co(H)$ and $g\in co(F/H)$ there exists
$h''\in co(H)$ such that $h\circ g = g\circ h''$. 
In particular there exist $h''_{1}, h''_{2} \in co(H)$
such that $h_{2}\circ g_{1} = g_{1}\circ h''_{2}$
and $h_{1}\circ g_{2} = g_{2}\circ h''_{1}$.
Then, 
\begin{equation}\label{equation1}
h_{1}\circ h_{2}\circ g_{1}\circ g_{2}(K)
=h_{1}\circ g_{1}\circ h''_{2}\circ g_{2}(K)
\subset h_{1}\circ g_{1}(K) = f_{1}(K).
\end{equation}
Because $co(H)$ and $co(F/H)$ are abelian (Lemma \ref{lem1}), we have 
$h_{1}\circ g_{1}\circ h''_{2}\circ g_{2}(K) 
= h_{2}\circ h_{1}\circ g_{1}\circ g_{2}$. Therefore,
\begin{equation}\label{equation2}
h_{1}\circ g_{1}\circ h''_{2}\circ g_{2}(K) 
= h_{2}\circ g_{2}\circ h''_{1}\circ g_{1}(K) 
\subset h_{2}\circ g_{2}(K) = f_{2}(K).
\end{equation}
Formulas (\ref{equation1}) and (\ref{equation2}) immediately 
show that $f_{1}(K)\cap f_{2}(K) \neq \emptyset$.
Again, an obvious induction proves that the family 
(\ref{family'}) has the finite intersection property. $\square$

\begin{lem} Let $F$ and $H$ be as above. Then, the family
\begin{equation}\label{family''}
\{f(K), f\in co(F) = co(H\circledS F/H)\},
\end{equation}      
of sets has the finite intersection property.\label{lem4} \end{lem} 
PROOF. Because the family $\{f(K), f\in \mathcal{R}\}$ has
the finite intersection property if and only if the family 
$\{f(K), f \in \mathcal{R}\cup \{I\}\}$ does, we can
add the identity mapping $I$ to the respective semi-group without loosing 
generality. By virtue of Lemma \ref{lem3} there exists a 
point $p$ common for all elements of the family (\ref{family'}).
In particular, $p \in h_{n}(K)$ and $p \in g_{n}(K)$ for all natural
$n$ and $h_{n}$, $g_{n}$ defined as follows\footnote{Here
$f^{n}$ denotes the $n$-fold composition 
$f\circ f \circ \ldots \circ f$  of $f$.}
\begin{displaymath}
h_{n} = \Big( \frac{1}{n}\Big(I + h + \ldots +h^{n-1}\Big)\Big)
\circ I 
= \frac{1}{n}\Big(I + h + \ldots +h^{n-1}\Big), \,\,\, h \in H,
\end{displaymath}
\begin{displaymath}
g_{n} = \Big( \frac{1}{n}\Big(I + g + \ldots +g^{n-1}\Big)\Big)
\circ I 
= \frac{1}{n}\Big(I + g + \ldots +g^{n-1}\Big), \,\,\, g \in F/H.
\end{displaymath} 
Let us take any $f \in F$. We have $f=h\circ g$, where $h\in H$ 
and $g \in F/H$, by assumption. It follows from what has already
been said that there exist points $x_{n}, y_{n} \in K$, such that
\begin{displaymath}
p = \frac{1}{n}\Big(x_{n} + h(x_{n}) + \ldots h^{n-1}(x_{n}))\Big),
\end{displaymath}
and
\begin{displaymath}
p = \frac{1}{n}\Big(y_{n} + g(y_{n}) + \ldots g^{n-1}(y_{n}))\Big).
\end{displaymath}
Therefore,
\begin{displaymath}
p - h(p)= \frac{1}{n}\Big(x_{n}  - g^{n}(x_{n}))\Big) 
\in \frac{1}{n}\Big(K-K \Big),
\end{displaymath}
\begin{displaymath}
p - g(p)= \frac{1}{n}\Big(y_{n}  - g^{n}(y_{n}))\Big) 
\in \frac{1}{n}\Big(K-K \Big).
\end{displaymath}
On the other hand if one takes any zero neighborhood $V$, then
$K-K \subset nV$ for the appropriately large $n$,
because $K$ is compact and thus linearly bounded.
In virtue of continuity of $h$ and $g$ we, therefore,
must have $p = h(p)$ and $p = g(p)$. Thus
\begin{equation}\label{equation3}
p = f(p) = h\circ g(p),
\end{equation} 
for all $f \in F$. 
Now, let $f$, $g \in co(F) = co(H\circledS F/H)$. It immediately
follows from what we have just proved that 
$p\in f(K) \cap g(K) \neq \emptyset$. The simple induction
shows now, that the family (\ref{family''}) has the finite 
intersection property. $\square$

Repeating inductively our Lemmas we can show that for any 
sequence $F_{1}, F_{2}, \ldots$ of semi-groups which 
we have mentioned above the family (\ref{family''})
with $F=F_{i}$ has the finite intersection property; 
or equivalently that (\ref{family''}) with 
$F\in \mathcal{F}$ has the finite intersection property.
The argument used in the proof of formula (\ref{equation3}),
shows that each semi-group $F\in\mathcal{F}$ possesses
a fixed point, i.e. we have

\begin{twr}[Generalized Markov-Kakutani theorem]
Let $K$ be a nonempty compact convex subset of linear topological
space $X$. Let $\mathcal{F}$ be the family of semi-groups
of continuous affine mappings which transform $K$ into itself, 
fulfilling Definition \ref{defin}. Then, for each $F\in \mathcal{F}$
there exist a point $p \in K$ such that $f(p) = p$, for all 
$f \in F$.\label{kakutani} \end{twr}  

Consider an abelian semi-group $F$ of bounded linear transformations
of a normed linear space $X$.
There is an interesting theorem in Functional Analysis,
which says that any $F$-invariant functional possesses
$F$-invariant extension which preserve norm, compare theorem 
5.24 of \cite{Rudin}. With the help of Theorem \ref{kakutani} 
we can prove that this theorem remains true for 
$F \in \mathcal{F}$ of Definition \ref{defin}. That is,

\begin{twr}[Generalized Hahn-Banach invariant theorem] Let $Y$ 
be subspace of a normed space $X$, let $g \in Y^{*}$ and
let $F \subset \mathcal{B}(X)$ fulfills\\
(a) $F \in \mathcal{F}$,\\
(b) $T(Y) \subset Y$, $\Vert T  \Vert \leq 1$, $T \in F$,\\
(c) $g\circ T = g$, for all $T\in F$.\\
Then, there exists functional $G \in X^{*}$ such that $G=g$
on $Y$, $\Vert G \Vert = \Vert g \Vert$ and $G\circ T = G$
for all $T \in F$. \end{twr} 
In the proof of this theorem one uses the set (one can assume 
$\Vert g \Vert = 1$)
\begin{displaymath}
K=\{\Lambda \in X^{*}: \Vert \Lambda \Vert \leq 1, \, \Lambda = g 
 \,\,\textrm{on} \,\, Y\},
\end{displaymath} 
as the the set $K$, which is not empty on account
of Hahn-Banach theorem. With the help of Banach-Alaoglu theorem 
one can show that $K$ is weakly-* compact. For each $T\in F$,
we define the affine mapping: $\Lambda \to \Lambda \circ T$,
which transforms $K$ into itself. It can be shown, compare
the proof of theorem 5.24 in \cite{Rudin}, that this mapping 
is weakly-* continuous. Theorem \ref{kakutani} implies now
existence of such $G \in K$ that $G\circ T = G$, for all $T\in F$.

\end{document}